\documentclass[12pt]{article}

\usepackage{amsmath,amssymb,amsfonts,setspace}
\usepackage{mathtools}

\def\BP{{\mathbb P}}

\def\bc{{\boldsymbol  C}}

\def\proof{\noindent{\bf Proof\ \ }}

\newcommand{\BZ}{\mathbb{Z}}

\newcommand{\eq}{\begin{equation}}
\newcommand{\en}{\end{equation}}

\newtheorem{theorem}{Theorem}

\newtheorem{lemma}{Lemma}

\newtheorem{corollary}{Corollary}

\def\endpf{\hfill $\Box$ \vskip0.5cm}
\def \proof{\noindent{\it Proof.\ }}

\begin{document}

\author{Alexander Gnedin and Dudley Stark}
\title{Random Permutations and Queues
}
\date{}

\maketitle

\begin{abstract}
\noindent
Given a growth rule which sequentially constructs random permutations of increasing degree, the stochastic process version of the rencontre problem asks what is the limiting proportion of time that the permutation has no fixed points (singleton cycles). We show that the discrete-time Chinese Restaurant Process (CRP) does not exhibit this limit. We then consider the related embedding of the CRP in continuous time and thereby show that it does have this and other limits of the time averages. By this embedding the cycle structure of the permutation can be represented as a tandem of infinite-server queues. We use this connection to show how results from the queuing theory can be interpreted in terms of the evolution of the cycle counts of permutations.

\end{abstract}

%\newpage

%\doublespacing

\section{Introduction} 
The { rencontre problem}  is one of  the oldest problems of probability theory.
%\cite{Ethier}.
If $n$ people exchange their hats at random, the probability that no
one receives their own hat is about $e^{-1}$.
% whenever  $n$ is not too small.
%equals $\sum_{j=1}^n (-1)^j/j!$, 
This classic result relates to the probability that the permutation of $[n]:=\{1,\ldots,n\}$ 
%obtained by splitting random permutation into cycles 
has no fixed points (singleton cycles). Such a permutation is called a derangement.
For permutations of increasing degree one can consider a dynamic version of  the rencontre problem. Given a
stochastic process that sequentially constructs permutations on $n=1,2,\ldots$ elements, 
how often a derangement will be observed?
% what is   the  limiting proportion of time 
%that the permutation has no fixed points? 
What are other path properties of the number of  fixed points seen as a process?
Questions of this sort 
can be asked also about  other functionals  related to the  cycle structure  of  the permutation.

%For instance, for the sequence of independent uniformly  distributed  random permutations, the limiting proportion  that the cycle partition has  sigletons about $e^{-1}$ 

A   growth rule which received considerable attention in the literature is  the following preferential attachment   algorithm known as the Chinese Restaurant Process (CRP) \cite{ABT, Feng, CSP}. 
Given a permutation of  $[n]$  which has been constructed at step $n$,
element $n+1$ is  either appended to the permutation as a new singleton cycle with probability proportional to $\theta$, or  inserted in random position within 
any  existing  cycle of size $m$ with probability proportional to $m$.
The  permutation obtained at step $n$  has the Ewens distribution, which is the uniform distribution on $n!$ permutations in the case $\theta=1$.
Properties of the Ewens distribution of fixed degree have been thoroughly studied. Still, there does not seem to be much work done on the dynamic properties connecting permutations
with variable $n$.

We will show that for the discrete-time CRP the  proportion of time when the permutation has no fixed points   does not converge.  An intuitive explanation for this 
is that the process counting fixed points slows down as $n$ increases, spending more and more steps  at the same level. 
To achieve convergence, we  employ a known  embedding of the CRP in  a continuous-time  birth process with immigration \cite{JT}.
In this realisation  the process counting fixed points
becomes identical with a  ${\rm M/M/}\infty$ queue, whose features have been intensely studied.
%  both in transient and stationary regimes.
Going deeper in this vein, the full process of cycle counts   behaves  as a series of such queues arranged in a
 tandem.
The principal point of the present note is that
this analogy opens  the way to 
 translate many results from the queueing theory 
 in terms of the evolutionary properties  of  permutations. There is vast literature 
 on infinite-server queues, so we do not attempt to survey the field exhaustively.
Much more our strategy is to collect and complement results allowing for transparent combinatorial interpretation,
with the primary focus on the small cycle counts of permutation.

\section{Background}

Let $\Pi^{(n)}$ be the  random permutation of $[n]$  at the $n$th step of the CRP.
The distribution of $\Pi^{(n)}$ is invariant under conjugations.
For  different degrees the permutations  are consistent, in the sense that $\Pi^{(m)}$, for  $m<n$,  can be derived from $\Pi^{(n)}$ by removing elements $m+1,\ldots,n$ from their cycles and  deleting empty cycles if necessary.

Let
$C_i^{(n)}$ be the
number of  cycles of size $i$ in $\Pi^{(n)}$. Thus  $C_1^{(n)}$ counts singletons (fixed points), $C_2^{(n)}$ doubletons, and so on.
The vector  of counts $\bc^{(n)}:=(C_1^{(n)},C_2^{(n)},\ldots)$ is a random  integer partition representing  the cycle structure of $\Pi^{(n)}$.
The generic value of $\bc^{(n)}$ is a vector ${\boldsymbol c}=(c_1, c_2,\ldots) \in \BZ_+^\infty$ satisfying
$\sum_{i=1}^\infty ic_i=n$
(hence $c_i=0$ for $i>n$). 
The law of  $\Pi^{(n)}$ is Ewens' distribution, which   assigns probability 
$\theta^c/(\theta)_n$ to 
permutation with  $c=\sum_{i=1}^n c_i$ cycles
(where $(\theta)_{n}=\theta(\theta+1)\cdots(\theta+n-1)$), so permutations of $[n]$ with the same number of cycles are equally likely.

%Given $\bc^{(n)}={\boldsymbol c}$ all  values of $\Pi^{(n)}$  with this cycle structure are equally likely.
%Partitions $(\Pi^{(n)}, ~n\geq1)$  are exchangeable and consistent. 
%The exchangeability means that all values of $\Pi^{(n)}$ 
%with given value $\bc^{(n)}={\boldsymbol c}$ are equally likely.
%The consistency means that $\Pi^{(n-1)}$  is obtainable fr%%om $\Pi^{(n)}$ by removing  the elements $n$ from its block (and deleting empty block if necessary).
%These properties imply that distribution of $\bc^{(n-1)}$ can be derived canonically from  distribution  $\bc^{(n)}$.

%thus the whole sequence $\Pi^{(1)},\cdots, \Pi^{(n)}$ is uniquely determined by $\Pi^{(n)}$.

The  process $(\bc^{(n)}, ~n\geq1)$  (we shall also use shorthand notation  $\bc^{(\cdot)}$) is a  nonhomogeneous Markov chain starting from the identity
 permutation for $n=1$ and
evolving in  $\BZ_+^\infty$  with transition probabilities
\begin{eqnarray}
\BP\big[\bc^{(n+1)}=(c_1+1,c_2,\ldots)|\bc^{(n)}={\boldsymbol c}\big]&=&
\frac{\theta}{\theta+n}, \label{trans1}\\
\BP\big[\bc^{(n+1)}=(c_1,\ldots,c_i-1,c_{i+1}+1,\ldots)|\bc^{(n)}={\boldsymbol c}\big]&=&
\frac{ic_i}{\theta+n},\quad c_i>0.\label{trans2}
\end{eqnarray}
Transition of the first type occurs when  $n+1$ starts   a new cycle, 
and of the second when the element is inserted in an existing cycle of $\Pi^{(n)}$.
The state distribution  is widely known as the Ewens sampling formula
\begin{equation}\label{ESF}
{\mathbb P}[\bc^{(n)}={\boldsymbol c}]=
\frac{n!}{(\theta)_{n}}   \prod_{i=1}^n 
\left(\frac{\theta}{i}   \right)^{c_i} \frac{1}{c_i!},~~~\sum_{i=1}^n i c_i=n.
\end{equation}
%where $(\theta)_{n}$ denotes the rising factorial.

The number of cycles in $\Pi^{(n)}$, denoted by
$$K^{(n)}:=\sum_{i=1}^n C_i^{(n)},$$
has probability generating function (p.g.f.)
\begin{equation}\label{KS}
{\mathbb E} z^{K_n}=\frac{(\theta z)_n}{(\theta)_n},
\end{equation}
which corresponds to the distribution of  the sum of $n$ independent Bernoulli variables with success probabilities $i/(i+\theta-1)$, $i=1,2,\ldots,n$.
As $n\to\infty$,
\begin{equation}\label{Kn}
K^{(n)}\sim \theta\log n~~{\rm a.s.},~~~\frac{K^{(n)}-\theta\log n}{\sqrt{\theta\log n}}\stackrel{d}{\to} {\rm N}(0,1).
\end{equation}
The counts $C_k^{(n)}$ for $k=1,2,\ldots,n$ are not independent because of the constraint in (\ref{ESF}).
Nevertheless, the small cycle counts of large permutation are almost independent: 

\begin{equation}\label{CZn}
\bc^{(n)}\stackrel{d}{\to}    (Z_1, Z_2,\ldots), ~~~n\to\infty,
\end{equation}
where the random variables $Z_k$ are independent with Poisson distribution
$Z_k\stackrel{d}{=} {\rm Poiss}(\theta/k)$. The convergence ({\ref{CZn}) holds with all moments.
See \cite{ABT} for detailed discussion including estimates of the convergence rate.

\section{Time averages in the CRP}

For $k\geq 1$ we denote
$\bc^{(n)}_k:=(C_1^{(n)},\ldots,C_k^{(n)})$  the truncated vector of the first $k$ cycle counts.
Easily  from (\ref{trans1})  and (\ref{trans2}),  $\bc^{(\cdot)}_k$ itself evolves as a Markov chain.
Each transition of  type (\ref{trans1})  or (\ref{trans2}) with $i\leq k$ triggers a jump of the truncated chain, while
 transitions of type (\ref{trans2}) with $i>k$ result in a loop, by which we mean same value $\bc^{(n+1)}_k=\bc^{(n)}_k$.
We assert that the average time 
 spent in any given state $(c_1,\ldots, c_k)$, that is
$\#\{m\leq n: \bc_k^{(m)}=(c_1,\ldots,c_k)\}/n$,
  does not have a limit as $n\to\infty$.
For the sake of simplicity of exposition we only consider the singleton count, the general case being completely analogous.

The process $(C_1^{(n)}, ~n\geq1)$
 is a nonhomogeneous Markov chain on ${\mathbb Z}_+$  which in  state $c$ has transition probabilities
\begin{eqnarray}\label{TraPro}
{\mathbb P} [C_1^{(n+1)}=j\,|\,C_1^{(n)}=c ]
=\begin{dcases*}
 \frac{\theta}{\theta+n}, ~~~{\rm for~} j=c+1, ~0\leq c\leq n,\\
 \frac{c}{\theta+n},  ~~~{\rm for~} j=c-1,~1\leq c\leq n,\\
 \frac{n-c}{\theta+n},  ~~~{\rm for~} j=c,~0\leq c\leq n.
\end{dcases*}
\end{eqnarray}
Loops occur in the event
   $C^{(n+1)}_1=C^{(n)}_1$.    
By time $n$ there will be about $K^{(n)}\sim \theta\log n$ upward moves caused by  starters, and about the same number of downward moves
caused by upgrading of singletons to doubletons.

We introduce here the time spent in state $c$,
$$T_n(c):= \#\{m\leq n: ~C_1^{(m)}=c\},$$
the time spent in state $c$, 
which equals
 the number of permutations among $\Pi^{(1)},\ldots,\Pi^{(n)}$ with exactly $c$ fixed points.
We will look at  the long-run  behaviour of the proportion $n^{-1} T_n(c)$ for fixed $c\geq 0$.
Specialising (\ref{CZn}) to singletons, we have  
\begin{equation}\label{eTheta}
{\mathbb P}[C^{(n)}_1
=i]\to  \frac{e^{-\theta} \theta^c}{c!}    , ~~n\to\infty.
\end{equation}
Taking the C{\'e}saro average in (\ref{eTheta}) results in the limit for the mean
${\mathbb E}[n^{-1} T_n(c)]\to  {e^{-\theta} \theta^c}/{c!}$,
which suggests that 
 the proportion itself obeys the law of large 
numbers.
But this intuition is wrong.

\begin{theorem}\label{CRP}
It holds that
$$
\liminf_{n\to\infty} n^{-1} T_n(c)    =0 ~~{\rm a.s.}~~{\rm and}~~ \limsup_{n\to\infty}  n^{-1} T_n(c) =1~~{\rm a.s.}
$$
\end{theorem}

\proof 
% Let $1=\tau_1<\tau_2<\dots$ be the sequence of times when the number of singleton changes, hence  $J_k:=C_1^{(\tau_k)}, k\geq1,$
% is the embedded jump chain.
Let $\nu_1<\nu_2<\ldots$ be the consequitive times when sojourns at $c$ start, that is $C_1^{(\nu_i-1)}\neq c$, and   $C_1^{(\nu_i)}=c$.
Choose an arbitrary integer $\gamma>1$, and consider the event
$A_k=\{\nu_{k+1}/\nu_k>\gamma\}.$
If $A_k$ occurs, the proportion $n^{-1}T_n(c)$ exceeds $1-1/\gamma$ for $n=\nu_{k+1}$.

Given $C_1^{(n)}=c$ and $m>n$, the probability that $C_1^{(j)}=c$ for all  $j=n,\ldots  ,m$ is 
\begin{eqnarray}\label{calc}
\prod_{j=n}^{m -1} \left(1- \frac{c+\theta}{j+\theta}\right)&=&
\exp\left(\sum_{j=n}^{ m -1}     \log \left( 1-\frac{c+\theta}{j+\theta} \right) \right)\nonumber\\
&=&
\exp\left(-(c+\theta)\sum_{j=n}^{m -1}  \frac{1}{j+\theta} +O\left(\frac{1}{n}\right)\right)\nonumber\\
&=&
\exp\left(-(c+\theta)\log\left(\frac{m}{n}\right)+O\left(\frac{1}{n}\right)\right).
\end{eqnarray}
Taking $m=\gamma n$, we obtain
$$
\prod_{j=n}^{\gamma n -1} \left(1- \frac{c+\theta}{j+\theta}\right)\to \gamma^{-(c+\theta)}\quad {\rm as~}n\to\infty.
$$
Using the strong Markov property to replace fixed $m$ by the random stopping time $\nu_k$, we conclude that ${\mathbb P}[A_k|{\Pi}_{\nu_{k}}]$ is bounded away from
$0$ as $k\to\infty$.
It follows that
$$\sum_{k=1}^\infty {\mathbb P}[A_k|{\Pi}_{\nu_{k}}]=\infty~~~{\rm a.s.}$$
Noting that $A_k$ is $\sigma({\Pi}_{\nu_{k+1}})$-measurable,
L{\'e}vy's conditional Borel-Cantelli lemma 
(\cite{Kallenberg}, p. 108)  applies to the sequence $(A_k, \sigma({\Pi}_{\nu_{k+1}}))$ and ensures that these events coincide:
$$\left\{\sum_{k=1}^\infty {\mathbb P}[A_k|{\Pi}_{\nu_{k}}]=\infty\right\}=\{A_k~{\rm i.o.}\}.$$
It follows that ${\mathbb P}[A_k~{\rm i.o.}]=1$.

Letting  $\gamma\to\infty$ we arrive at
$\limsup_{n\to\infty} n^{-1} T_n(c) =1\,{\rm~ a.s.}$
for every $c$.
But then since $n^{-1} T_n(c) + n^{-1} T_n(d)\leq1$   for $d\neq c$
we also have
$\liminf_{n\to\infty} n^{-1} T_n(c)=0 {\rm ~~a.s.}$
\endpf

\vskip0.2cm
We leave to the reader checking that the variance of $n^{-1}T_n(c)$ does not vanish asymptotically, hence convergence  in probability also fails.
The source of irregularity  is the infinite expectation of sojourn times.
This phenomenon is akin to the null recurrence of  time-homogeneous Markov chains like the symmetric random walk.

\section{The embedded random walk} \label{RW}

The convergence of time averages can be achieved by discarding the loops, that is letting the clock tick
only at times of nontrivial moves.

To explore this thread, 
let $(J_i,~i\geq0)$ be a  nearest-neighbour random walk on ${\mathbb Z}_+=\{0,1,2,\ldots\}$ which moves  $\pm1$ with time-independent
probabilities
\begin{equation}\label{pq}
p_c= \frac{\rho}{c+\rho},~~~~~q_c=\frac{c}{c+\rho},
\end{equation}
respectively,
where $\rho>0$ is a parameter.
The random walk is reversible and, checking the detailed balance equations, it is  seen that 
it has the unique stationary  distribution
\begin{equation}\label{alpha-c}
\alpha_c= \frac{e^{-\rho}(\rho+c)\rho^{c-1}}{2\,\, c!}, ~~c\in{\mathbb Z}_+,
\end{equation}
(which can be decomposed as  equi-weighted   mixture of  Poiss$(\rho)$ distribution on $\{0,1,\ldots\}$ and the shifted Poiss$(\rho)$ distribution on $\{1,2,\ldots\}$).
Appealing to the ergodic theorem, we conclude that  the proportion of time that $J_\cdot$ spends at $0$ converges to $\alpha_0=e^{-\rho}/2$.
Therefore 
the mean time between two consequitive visits to $0$ is $1/\alpha_0= 2 e^\rho$, and
the mean number of upward moves between the visits  is $e^\rho$.

Define an excursion of $J_\cdot$  above level $c$ to be a segment of the path that starts from
$c+1$ and terminates by hitting $c$.
Assuming $J_0=c+1$, the first passage time

\begin{equation}\label{mean-kappa}
\kappa_c:= \inf\{k\geq 0: J_k=c\}
\end{equation}
is the length, and
$$H_c:=\sup\{J_k: k\leq\kappa_c\}-c$$
is  the height of excursion above $c$.

The random walk starts anew by each  visit at $c$, hence the elementary renewal theorem ensures that  the number of   excursions above $c$ 
completed within $k$ steps is
 asymptotic to
$$\alpha_c\frac{\rho}{\rho+c}\, k,~~~k\to\infty, $$
which entails that the limit proportion of time spent above $c$ is
$$\alpha_c\frac{\rho}{\rho+c}\, {\mathbb E}[\kappa_c],$$
which must also be equal to
$\sum_{j=c+1}^\infty \alpha_j$
by ergodicity.
Recalling (\ref{alpha-c}), an easy calculation gives
$${\mathbb E}[\kappa_c]= 1+\frac{2c!}{\rho^c} \sum_{j=c+1}^\infty \frac{\rho^j}{j!}.$$
In particular,
$${\mathbb E}[\kappa_0]=2e^\rho-1.$$

To determine  the variance of  excursion length  we adopt a formula of Harris (\cite{Harris}, Equation (5.9)). To that end,
express the  product of odds involved in the cited result in terms of Poiss$(\rho)$ probabilities as
\begin{equation*}
\pi_r:= e^{-\rho} \frac{\rho^{r}}{r!},~~~{\rm so~~}\prod_{j=1}^{r-1}\frac{q_j}{p_j}= (e^\rho\pi_{r-1})^{-1},
\end{equation*}
to obtain
\begin{equation}\label{Var-D}
{\rm Var}\, [\kappa_0]= 4e^\rho(2\rho-e^\rho+1) +8e^\rho\sum_{r=0}^\infty \frac{1}{\pi_r}\left(\sum_{j=r+1}^\infty\pi_j\right)^2.
\end{equation}

Harris (\cite{Harris}, Theorem 2b) also solved the generalised gambler's ruin problem:
for fixed integers  $0\leq \ell<s\leq u$,
 if the random walk starts at $s$ it will reach $u$ before visiting $\ell$ with probability
\begin{equation}\label{Harris}
{\mathbb P}[J_\cdot {\rm ~reaches~} u {\rm ~before~}\ell|J_0=s]=
\frac{\cfrac{1}{\pi_\ell}+\cfrac{1}{\pi_{\ell+1}}  +\cdots+ \cfrac{1}{\pi_{s-1}}}     {\cfrac{1}{\pi_\ell}+\cfrac{1}{\pi_{\ell+1}}  +\cdots+ \cfrac{1}{\pi_{u-1}}} 
\end{equation}
Choosing $\ell=c, s=c+1, u=h+1$ this gives the distribution of $H_c$, most conveniently expressed in the form of the upper tail probabilities
\begin{equation}\label{height}
{\mathbb P}[H_c\geq h+1]=\frac{1}{\sum\limits_{r=0}^h \cfrac{(c+1)_r}{\rho^r}}, ~~~h\geq0.
\end{equation}

The tails are lighter than geometric,
$$\lim_{h\to\infty}\frac{ {\mathbb P}[H_0\geq h+1]}{{\mathbb P}[H_0\geq h]}= 0,$$
which suggests that the maximum of random walk satisfies a law of large numbers for discrete random variables  \cite{Anderson}. This was indeed shown by 
Park et al \cite{Park}. To state their exceptionally precise result, let $\beta_m$ be the time of  the $m$th visit of the random walk to $0$. Then
\begin{equation}\label{Park}
\lim_{m\to\infty} {\mathbb P}\,[\max_{0\leq i\leq\beta_m} J_i\in \{ I_m, I_m+1\}]=1,
\end{equation}
where 
\begin{equation}\label{Iam}
I_m:=\left\lfloor \frac{\log m-\frac{1}{2}\log\log m-\frac{1}{2}\log 2\pi}{\log\log m-1-\log\rho}   +\frac{1}{2}      \right\rfloor.
\end{equation}

Using $\rho=1$ in (\ref{height}), we obtain numerical values for the excursion height moments
\begin{eqnarray*}
{\mathbb E} [H_0]=\sum_{h=0}^\infty \frac{1}{0!+1!+2!+\cdots+h!}=1.887\ldots, \\
{\rm Var} [H_0]=\sum_{h=0}^\infty \frac{2 h+1}{0!+1!+2!+\cdots+h!}-({\mathbb E}[H_0])^2=1.242\ldots
\end{eqnarray*}
and the number of upward moves of an excursion has expectation and variance
$${\mathbb E}\left [\frac{\kappa_0-1}{2}\right]=1.718\ldots,~~~ {\rm Var}\left[\frac{\kappa_0-1}{2}\right]=7.930\ldots $$

Translating the results above to the CRP and taking our parameter $\rho=\theta$, the random walk $J_\cdot$ is the embedded jump chain for $C_1^{(\cdot)}$, 
hence we can make conclusions about the number of fixed points in the CRP.
Excursions above $c=0$ correspond to the fluctuation in the number of fixed points in the period between two consequitive derangements.
The asymptotic proportion of derangements within the number of nontrivial moves of $C_1^{(\cdot)}$ is $e^{-\theta}/2$.
This  does not match with  the value $e^{-\theta}$  that could be anticipated from (\ref{eTheta}).
Starting from $C_1^{(n)}=c+1, ~c\geq0,$ the variable $H_c+c$ is  the maximum number of singletons observed until their number falls to $c$, 
and $(\kappa_c-1)/2$ is the  number of new  cycles produced by the CRP within this period.
In particular, applying the above findings  to the case $\theta=1$ of uniformly chosen permutations, we see that a permutation with initially one singleton (for instance, $\Pi^{(1)}$) will have on the average 
$1.778$ cycles at the first time when it becomes derangement, and the expected maximum number of singletons observed by this time is about $1.887$.

For the first $n$ permutations $\Pi^{(1)},\ldots,\Pi^{(n)}$,
the number of fixed points $C^{(\cdot)}_1$ will change the value about $2\theta\log n$ times, when singletons are formed and when they progress to doubletons.
This implies a strong law of large numbers  on the log scale
$$\#\{j\leq n: \Pi^{(j)} ~{\rm is~a~derangement~and~}\Pi^{(j-1)}~{\rm is~not}\}\sim (2\theta\log n)\alpha_0= e^{-\theta}\theta \log n~~~{\rm a.s.}$$
Thus  $e^{-\theta}$ appears to be the asymptotic proportion of singletons that enter
the permutation when it is a derangement, relative to the number of all cycles.
 For the maximum number of fixed points we have from (\ref{Kn}) and  (\ref{Park}) 
$$\lim_{n\to\infty} {\mathbb P}\,[ \max_{1\leq j\leq n} C^{(j)}_1\in \{I_m,I_m+1\}]=1,$$
where $I_m$ is given by (\ref{Iam}) with a possible adjustment of $\pm 1$,
with $\rho=\theta$ for 
$m=\lfloor\theta\log n\rfloor$. 
The fact that approximating the true value of $m$ by  $\lfloor\theta e^{-\theta}\log n\rfloor$ only changes the expression inside $\lfloor\cdot\rfloor$ in  (\ref{Iam}) by $o(1)$ accounts for the possible adjustment to $I_m$.

It is tempting to similarly link $C_k^{(\cdot)}$ to the random walk $J_\cdot$ 
with parameter $\rho=\theta/k$. For  $k>1$,  this comparison does not work literally, because
the count of  $k$-cycles is not a Markov chain. The connection becomes valid  asymptotically, in the sense that
for $n_0\to\infty$ the loop-free path of $(C_k^{(n)},~n\geq n_0)$
%(seen as a functional of Markovian $\bc^{(\cdot)}_k$)
conditioned on $C_k^{(n_0)}=c$ converges in distribution to $J_\cdot$ with the initial state $J_0=c$.
This will follow from the embedding in the next section. To gain some intuition, observe that the transition of $\bc_{k-1}^{(\cdot)}$  causing an upward move of $C_k^{(\cdot)}$ 
has  probability
$${\mathbb P}[\bc_{k-1}^{(n+1)}-\bc_{k-1}^{(n)}=(0,\ldots,0,-1)\,|\,\Pi^{(n)}]= \frac{(k-1)C_{k-1}^{(n)}}{n+\theta}.$$
For large $n$ the distribution of  $C_{k-1}^{(n)}$ is approximately Poisson with mean $\theta/(k-1)$, hence the unconditional probability of the said upward  move is about $\theta/(n+\theta)$,
to be compared with probability $kc/(\theta+n)$ of the downward move (which does not depend on the first $k-1$ counts); thus for  large times given $C_k^{(\cdot)}$ has a move 
it is $+1$
with probability about $(\theta/k)/(c+\theta/k)$, in agreement with (\ref{pq}).

\section{Embedding in continuous time}\label{embed}

To ensure the convergence of  time averages of occupation times, the temporal scale of the CRP should be changed so
that the degree of permutation grows  about exponentially, and hence the number of cycles grows about linearly.
An elegant way to do this is to embed 
the permutation-valued process in continuous time. The embedding idea
originated in  \cite{AT,  JT,  Tav} and is nicely presented in \cite{Feng}.

Consider  a permutation-valued process $(\Pi(t),~t\geq0)$ which starts with the empty permutation (of degree $0$) and then
evolves according to  this rule: 
 given permutation $\Pi(t)$ of $[n]$, element $n+1$ starts a new cycle with probability rate 
$\theta$ and is inserted in random position of any existing cycle of size $m$ at rate $m$.
It is obvious from this description that the associated discrete-time jump chain  is the CRP $(\Pi^{(n)},~n\geq0)$.

%Partitions $\Pi(t)$ are still consistent in a weak sense. Fix $s<t$. Let $p=(e^s-1)/(e^t-1)$ and consider a thinning procedure which retains each element of $\Pi(t)$ with probability $p$.
%Enumerate the retained elements by integers in increasing order: the resulting partition has same distribution as $\Pi(s)$.

Let
  ${\boldsymbol C}(t)=(C_1(t), C_2(t),\ldots)$ denote
the    vector counting singletons, doubletons, and so on. 
Clearly,
  $({\boldsymbol C}(t), t\geq0)$ is a Markov process 
on  $\{{\boldsymbol c}\in{\mathbb Z}_+^\infty:~\sum_i c_i<\infty\}$ 
with time-independent transition rates
\begin{eqnarray*}
\theta ~~{\rm for~~}  (c_1, c_2,\ldots)&\to& (c_1+1, c_2,\ldots),\\
i c_i  ~~{\rm for~~}  (c_1, c_2,\ldots)&\to& (c_1, \ldots, c_i-1, c_{i+1}+1,\ldots), ~ i\geq2,
\end{eqnarray*}
and (by definition) the initial value ${\boldsymbol C}(0)=(0,0,\ldots)$.

Define $K(t)$ and $N(t)$ to be the number of cycles and  the degree of permutation 
$\Pi(t)$, respectively. Thus,
$$K(t):=\sum_{i=1}^\infty C_i(t)~~~{\rm and}~~~N(t):=\sum_{i=1}^\infty i\, C_i(t).$$
The number of cycles  evolves according to a Poisson process of rate $\theta$.
The degree of $\Pi(t)$ follows  a linear birth process with immigration, where  the immigration rate is constant $\theta$ and the  birth rate per capita is $1$
(so that given $N(t)=n$, the birth rate is $n$).
Sometimes the term Pascal process is used for such a process, because
 the conditional distribution of  the increment $N(t)-N(s)$ given $N(s)=n$ is negative binomial (i.e.~Pascal)
${\rm NB}(\theta+n, 1-e^{-(t-s)})$, for  $t>s\geq0$.
See the texts \cite{Asmussen, Feng, Resnick} for properties of the birth-death processes.

%\item Tandem of  A network consists of  service stations $1,2,\ldots$, each with infinitely many servers and {\rm Exponential}$(k)$ service time at station $k$.
%Customers enter station 1 at rate $\theta$, then after getting served proceed to station 2, etc. Component $C_i(t)$ is the number of customers at station $i$.
%Parameter $\theta$ appears here as the workload of the queues.
%\end{itemize}

%\begin{itemize}
%\item A linear birth process with immigration. 
%The migrants enter the population at rate $\theta$ and initiate families. Every individual in the existing population produces offsping at unit rate.
%Component $C_i(t)$ is the number of families of size $i$.

We let ${\boldsymbol  C}_k(t):=(C_1(t),\ldots, C_k(t))$ denote a truncated vector of  cycle counts. 
Similarly to the discrete-time CRP,
each ${\boldsymbol  C}_k(\cdot)$
 is itself a time-homogeneous Markov process.

The next product-form  result is known in much larger generality in the theory of networks and population processes \cite{Kelly}.  
The textbook proofs for the transient (pre-limit) state all employ Kolmogorov's equation.
The following elegant elementary proof for the special case in focus was  outlined in  \cite{AT}, Exercise 10.7.

\begin{theorem}\label{IndBl} The random variables $C_1(t), C_2(t),\ldots$ are independent, with  distribution
\begin{equation}\label{transient}
C_k(t)\stackrel{d}{=} {\rm Poisson}\left(  \frac{\theta(1-e^{-t})^k  }{k} \right).
\end{equation}
Therefore, as $t\to\infty$
\begin{equation}\label{steady}
{\boldsymbol C}(t)\stackrel{d}{\to} (Z_1, Z_2,\ldots),
\end{equation}
where $Z_i\stackrel{d}{=}{\rm Poisson}(\theta/i)$ are independent.

\end{theorem}

\proof  A singleton needs time $\xi_1$ to become a doubleton, then time $\xi_2/2$ to become a tripleton, etc.,
where $\xi_1,\xi_2,\ldots$ are independent unit exponential random variables.
The jomp times for different cycles are independent.
Hence by the theorem on marked Poisson processes (see \cite{KingmanPP}, Ch. 5),
the $C_k(t)$'s for $k=1,2,\ldots$ are independent and Poisson distributed. It remains to compute the means.

A singleton grows to a cycle of size at least  $k+1$ within time  $s$  with probability 
\begin{equation}\label{OrderSt}
{\mathbb P}\left[\sum_{i=1}^k {\xi_i}/{i} <s\right]={\mathbb P}[\max(\xi_1,\ldots,\xi_k)<s]=(1-e^{-s})^k,
\end{equation}
where the first identity  is R{\'e}nyi's representation of the exponential order statistics.
Hence, the probability that this is a $k$-cycle  is 
$$(1-e^{-s})^{k-1}- (1-e^{-s})^k= e^{-s} (1-e^{-s})^{k-1}.$$
Thus, the mean number of  $k$-cycles at time $t$ is 
$${\mathbb E}\, [C_k(t)]=\int_0^t       e^{-s} (1-e^{-s})^{k-1}         \theta\,{\rm d} s=\theta  k^{-1}(1-e^{-t})^k,$$
as wanted.

The convergence of ${\boldsymbol C}(t)$ in distribution follows from the first part of the statement.
Note that this convergence
is understood relative to the discrete product topology on ${\mathbb Z}_+^\infty$ and amounts to the weak convergence of truncated processes
${\boldsymbol C}_k(t)$.
\endpf

The Markov chain $({\boldsymbol C}_k(t), ~t\geq0)$ is positive recurrent, hence application of the ergodic theorem ensures  existence of  the time averages.
Let
$$h_k:=\sum_{i=1}^k \frac{1}{i}$$
denote the $k$th harmonic number.

\begin{corollary} The proportion of time spent by  the process $({\boldsymbol C}_k(t), ~t\geq0)$
in state $(c_1,\ldots,c_k)$ converges to 
$${\mathbb P}[Z_1=c_1,\ldots,Z_k=c_k]= e^{-{\theta}h_k } \prod_{i=1}^k \frac{ \theta^{c_i}}{i^{c_i} c_i!} \,.$$
\end{corollary}

\noindent
In particular, the average time when  $\Pi(\cdot)$ is  a derangement approaches  $e^{-\theta}$.

We may tag a cycle by its minimal element. The growth  of a cycle can be thought of as
passing through phases $S_1, S_2,\ldots$ of being a singleton, doubleton, etc. 
 The sojourn periods  across different cycles and phases are independent.
The input flow into  $S_1$ is a Poisson process of rate $\theta$, and the time spent in $S_k$  
has exponential distribution with parameter $k$. By Theorem \ref{IndBl} the flow from $S_k$ to $S_{k+1}$ is nonhomogeneous Poisson with rate $\theta(1-e^{-t})^k$,
hence converging to homogeneous flow with rate $\theta$.

Such a process,
with general  sojourn rates  $ \mu_1,\mu_2,\ldots$, models a network 
 of ${\rm M/M/}\infty$, infinite-server queues  connected in a tandem  \cite{Boxma, GS, KY, P}.
In the literature the tandem is often considered as open
 network  with finitely many phases $S_1, \ldots, S_k$,  where the task
 departs upon passing through $S_k$.

The process $({\boldsymbol C}_k(t),~t\geq 0)$  is stationary if  it starts  with the product Poisson distribution  ${\boldsymbol C}_k(0)\stackrel{d}{=} (Z_1,\ldots,Z_k)$ as in (\ref{steady}).
Following the established terminology we shall call the stationary process steady state, as opposed to the transient regime with the pre-limit 
state distribution (\ref{transient}). 
In the steady state  the flow from $S_i$ to $S_{i+1}$ is Poisson with rate $\theta$, so each $C_i(\cdot)$ behaves like a single stationary ${\rm M/M/}\infty$ queue;
this is an instance of the seminal Burke's theorem.
We stress that the steady state does not  describe permutations of finite degree, but rather captures asymptotic  features of  small cycle counts of  $\Pi(\cdot)$ at large times.

\section{Pascalisation and big cycles}

The discrete- and continuous time models are related via
$$(\Pi(t),~t\geq0)\stackrel{d}{=} (\Pi^{(N(t))}, t\geq0),$$
where the permutation in the right-hand side is constructed from two independent  ingredients: CRP $(\Pi^{(n)},~n\geq0)$ and a Pascal process $(N(t), ~t\geq0)$.
For this kind of randomisation we  propose the term {\it pascalisation}, by analogy with the established concept of poissonisation (sampling $n$ from the Poisson distribution).
These methods are most useful in the situations where they produce exact independence instead of 
 the asymptotic independence   in fixed-$n$ combinatorial models.
In the context of cycle structure  the pascalisation was used already in \cite{Shepp} for the case $\theta=1$ (where the mixing distribution is  geometric), 
in particular to prove the convergence (\ref{CZn}). However, for $\theta\neq 1$ the method does not seem to have been given due attention.
See \cite{BO} for pascalisation of another interesting distribution on integer partitions.

To illustrate,   let $g_n(z)$ be the p.g.f. of   $K^{(n)}$. Connecting $g_n(z)$ to the Poisson p.g.f. of $K(t)$ produces
$$ e^{\theta t(z-1)}=\sum_{n=0}^\infty {\mathbb P}[N(t)=n]\,g_n(z)=   \sum_{n=0}^\infty\frac{(\theta)_n}{n!}e^{-\theta t}(1-e^{-t})^n g_n(z).$$ 
Expanding the left-hand side in powers of $(1-e^{-t})$ and equating the coefficients yields
$g_n(z)= {(\theta z)_n}/{(\theta)_n}$, which gives altervative proof of (\ref{KS}). For both discrete and continuous models, the number of cycles is about normally distributed  for large times.

A more complex functional is the maximal  size of a cycle. For $\Pi^{(n)}$ this has a  sophisticated  limit distribution (\cite{Feng}, Theorem 2.5), but for $\Pi(t)$ 
the things are rather straightforward.
Let $M(t):=\max\{i: C_i(t)>0\}$ be the maximal  size of a cycle present in $\Pi(t)$. From (\ref{transient}),
$${\mathbb P}[M(t)\leq m]=\exp\left( -\sum_{k=m+1}^\infty \frac{\theta(1-e^{-t})^k}{k}   \right).$$
To find the limit law of $M(t)$ consider a Poisson random measure (PRM) ${\cal P}_t$ which charges 
point $e^{-t}k$ with mass $C_k(t)$, $k=1,2,\ldots$ ans let $\cal P$ be another PRM on $(0,\infty)$ with mean measure $\lambda({\rm d}x)=\theta e^{-x} {\rm d}x/x$.
We have  $\lambda (x,\infty)=\theta E_1(x)$, where 
$$E_1(x)=\int_x^\infty \frac{e^{-s}}{s}{\rm d}s$$ 
is the exponential integral function.

\begin{theorem} For $t\to\infty$, the PRM ${\cal P}_t$ converges weakly to $\cal P$.
\end{theorem}

\proof It is sufficient to show that the mean measure $\lambda_t$ of ${\cal P}_t$ satisfies
$$\lim_{t\to\infty} \lambda_t(x,\infty)=\lambda(x,\infty)$$
 for each $x>0$. We have
\begin{eqnarray*}
\lambda_t(me^{-t},\infty)=\theta \sum_{k=m+1}^\infty  \frac{(1-e^{-t})^k}{k}&=& \theta \sum_{k=m+1}^\infty \int_{e^{-t}}^1 (1-y)^{k-1}{\rm d}y=\\
\theta \int_{e^{-t}}^1 \frac{(1-y)^m}{y} {\rm d}y&=&\theta \int_{1}^{e^t} \frac{(1-se^{-t})^m}{s}{\rm d}s.
\end{eqnarray*}
Setting $m=\lfloor xe^t\rfloor $ we obtain, by the monotone convergence, 
\begin{eqnarray*}
\lambda_t(x,\infty)=\theta \int_1^{e^t} \frac{(1-se^{-t})^{xe^t}}{s}{\rm d}s +o(1)\to 
\theta \int_1^{\infty} \frac{e^{-xs}}{s}{\rm d}s= \theta E_1(x),
\end{eqnarray*}
and the conclusion follows.
\endpf
The result implies that   $e^{-t}M(t)$ converges in distribution to the largest point of $\cal P$, whence
$$\lim_{t\to\infty} {\mathbb P}[e^{-t} M(t)\leq x]= \exp\left(-\theta E_1(x)\right).$$
Moreover,  the whole scaled decreasing sequence of the cycle lengths of $\Pi(t)$ converges in distribution to the infinite sequence of points of $\cal P$ listed in decreasing order.
The analogous limit for $\Pi^{(n)}$, the Poisson-Dirichlet distribution,  can be obtained by normalising the points of $\cal P$ by their sum (\cite{Feng}, Theorem 2.2).

Better tractable limit laws  for  large cycles in $\Pi^{(n)}$ or $\Pi(t)$ appear if the cycles are listed 
in the age order, that is  by increase of the singletons. For instance, the size of the oldest cycle of $\Pi(t)$ (containing element $1$) is asymptotic to $e^{t-\eta/\theta} \xi $
with independent unit exponential $\xi$ and $\eta$.
See \cite{ ABT-PD, ABT, Feng} for various representations  of the multivariate limit.

\section{Excursions of a cycle count process}\label{one}

The   ${\rm M/M/}\infty$ queue occupancy process   is a Markov chain  $(X(t), ~t\geq0)$ on ${\mathbb Z}_+$ which from state $c$ jumps by $\pm1$ with rates $\theta$ and $\mu c$, respectively.
This is sometimes called a linear immigration-death process \cite{Asmussen}. We shall follow the intuitive terminology of queueing theory, calling 
$\theta$ the arrival rate
of the input Poisson process $K(\cdot)$, $\mu$ the service rate (the departure rate per task),
and $\rho:=\theta/\mu$ the average workload.
The transient state distribution is Poisson with parameter depending on $t$, that is
$${\mathbb P}[X(t)=c\,|\,X(0)=0]=\pi_c(t),~~~{\rm where~~~}\pi_c(t)= \exp\{-\rho(1-e^{-t/\mu}\}\, \frac{[(\rho(1-e^{-t/\mu})]^c}{c!},~~~c\geq0,$$
and in the steady state the distribution is ${\rm Poiss}(\rho)$.
The embedded jump chain for $X(\cdot)$ is the random walk $J_\cdot$  with parameter $\rho$ as in Section \ref{RW}.

In terms of the permutation-valued process, $X(\cdot)$ could be $C_1(\cdot)$, or $C_k(\cdot)$ for $k> 1$ with Poisson inflow resulting  from the output  of    $C_{k-1}(\cdot)$ (or $\bc_{k-1}(\cdot)$) in the steady state. Then the parameters are $\theta, \mu=k, \rho=\theta/k$.

For $c\geq0$. we define  excursion above $c$ to be a segment of the path 
that starts at $c+1$ and terminates by the first passage of level $c$. 
The case $c=0$ is referred to as the busy period, and for the general $c$ the excursion is called the congestion period above the level.
Excursions below $c\geq 1$ are defined analogously but will not be touched here (see \cite{R} or \cite{Mandjes} on intercongestion periods).

A visit to $c$ is followed by an excursion above $c$ if the next state is $c+1$.
%(this is automatic for $c=0$). 
In the long run,  the  mean rate  of the point process of jumps $c\to c+1$ is about $\pi_c\lambda$,
thus by renewal theory the number of  excursions above $c$ completed by  time $t$ is asymptotic to $t/(\pi_c\lambda)= te^\rho c!/(\lambda \rho^c)$ as $t\to\infty$.

The functionals characterising  the excursion include  
%the duration, the overflow,  the height and  the number of new arrivals, which are, respectively
\begin{eqnarray*}
{\rm the~ duration~~~~}D_c&=&\inf\{t: X(t)=c\},\\
{\rm the~ height~ above~}c~~~H_c&=&\sup\{X(t)-c: 0\leq t\leq D_c\},\\
{\rm the~ overflow~~~~}A_c&=& \int_0^{D_c} (X(t)-c){\rm d}t,\\
{\rm the ~number~ of ~new~arrivals ~~~}\Delta_c&=&(\kappa_c-1)/2,
\end{eqnarray*}
where  we write definitions as if the excursion started   at time $0$ with $X(0)=c+1$. 
The variables $H_c$ and $\kappa_c$ are functionals of the embedded random walk $J_\cdot$ 
and have the same meaning as in Section \ref{RW}. In the rest of this section we put together and complement 
properties  of these variables found in the literature.

Some relations among the moments follow by the optional sampling theorem applied to 
$\Delta_c=K(D_c)$ and
the martingale $K(\cdot)-\theta \,\cdot$ where $K(\cdot)$ is the Poisson arrival process:
\begin{equation}\label{Wald}
{\mathbb E}\left[\Delta_c-\theta D_c\right]=0,~~~{\mathbb E}[\Delta_c-\theta D_c]^2=  {\mathbb E}[\Delta_c].
\end{equation}
By arguments from the renewal theory, 
\begin{equation}\label{RTA}
{\mathbb E} [D_c]=\frac{1}{\theta\pi_c}\sum_{j=c+1}^\infty \pi_j,~~~{\mathbb E} [A_c]=\frac{1}{\theta\pi_c}\sum_{j=c+1}^\infty \pi_j (j-c),
\end{equation}
and from (\ref{Wald}) 
\begin{equation}\label{Edelta}
{\mathbb E}[\Delta_c]=\frac{1}{\pi_c}\sum_{j=c+1}^\infty \pi_j,
\end{equation}
where $\pi_j$ are the Poisson$(\rho)$ probabilities.
The variance of $\Delta_c$ for $c=0$ is obvious from the connection with $\kappa_c$ and (\ref{Var-D}), and for $c>0$ can be also derived from (\ref{Harris});
from (\ref{Wald}) one can compute then the covariance between $D_c$ and $\Delta_c$.

Formulas  for the Laplace transforms of these statistics have been obtained in terms of the integrals
$$I_c(\alpha,\beta)=\int_0^1 u^{c} (1-u)^{\alpha-1}  e^{-\beta  u}{\rm d}u,$$
which in turn can be expressed through Kummer's confluent hypergeometric function 
$$M(a,b,z):=\sum_{i=0}^\infty \frac{(a)_i}{(b)_i}\frac{z^i}{i!}$$
as
$$I_c(\alpha,\beta)= e^{-\beta} \frac{\Gamma(c+1)\Gamma(\alpha)}{\Gamma(c+\alpha+1)}  M(\alpha,\alpha+c+1,\beta).$$
The appearance of these functions here is quite natural, since the Laplace transform of the transient state probability $\pi_c(t)$ is
$$\int_0^\infty \pi_c(t)e^{-z t}{\rm d}t= \frac{\rho^c\mu}{c!}I_c(\mu z,\rho),$$
as one can easily calculate by
\begin{eqnarray*}
\int_0^\infty \pi_c(t)e^{-z t}{\rm d}t&=&
\int_0^\infty exp(-\rho(1-e^{-t/\mu}) \frac{[\rho(1-e^{-t/\mu}) ]^c}{c!}e^{-zt}dt\\
&=&\int_0^1e^{-\rho u}\frac{[\rho u]^c}{c!} (1-u)^{\mu z} \mu(1-u)^{-1}du\\
&=&\frac{\mu\rho^c}{c!} \int_0^1e^{-\rho u}u^c (1-u)^{\mu z-1} du\\
&=&\frac{\mu\rho^c}{c!} I_c(\mu z,\rho),
\end{eqnarray*}
where $u=1-e^{-t/\mu}$, $du=\frac{1}{\mu}e^{-t/\mu}dt\Rightarrow dt=\mu(1-u)^{-1}du$
and $e^{-zt}=(1-u)^{\mu z}$.
Concretely, Guillemin and Simonian \cite{GS} showed that 
\begin{equation}\label{GS}
{\mathbb E}[\exp(-z D_c)]=\frac{I_{c+1}(z/\mu,\rho)}{I_{c}(z/\mu,\rho)}.
%=\frac{c+1}{c+1+z/\mu}\,\,\frac{M(z/\mu,c+z/\mu+2,\rho)}{M(z/\mu,c+z/\mu+1,\rho)}    ,
\end{equation}

Preater \cite{P}  derived the joint Laplace transform
\begin{equation}\label{Preater}
{\mathbb E}[\exp(-x D_0-y \Delta_0- z A_0)]=\frac{\mu}{z+\mu} \,\frac{I_{c+1}(a-b,b)}{I_{c}(a-b, b)},
\end{equation}
where 
$$a=\frac{x+\theta}{z+\mu},  ~~~b=\frac{\theta\mu e^{-y}}{(z+\mu)^2}$$
(this result is cited as Equation (20) in \cite{Mandjes}), and in  \cite{P0} obtained a continued fraction formula  for the joint Laplace transform of $D_c$ and $A_c$.
Preater's approach \cite{P0, P} to continued fractions expansions relies on the fact that
an excursion above $c$ decomposes in  a sojourn  at $c+1$ of
${\rm Exp}(\theta+(c+1)\mu)$-length and some ${\rm Geom}((c+1)/(\rho+(c+1)))$ number of
path segments, each comprised of excursion above $c+1$ and sojourn at $c$, with all ingredients  being independent. 

% in terms of continued fractions \cite{GP, P0, P, Mandjes} or

Roijers et al \cite{Mandjes} notice that  obtaining
  higher moments by
differentiating the Laplace transforms
is not straightforward due to the implicit nature of functions 
involved.
They  derived recursions in $c$ using the said decomposition of the excursion 
above $c$, thus eventually reducing to the case $c=0$.
For  the second moment of the duration they obtain a series expansion
 (\cite{Mandjes}, Equation (24))
\begin{equation}\label{D2}
{\mathbb E}[D_0^2]= 
\frac{2e^{2\rho}}{\theta\mu} \sum_{j=1}^\infty \frac{\pi_j}{j},
\end{equation}
which can be written as
$${\mathbb E}[D_0^2]=\frac{2e^\rho}{\theta\mu}\int_0^\rho \frac{e^s-1}{s} {\rm d}s.$$
Lizgin and Rudenko \cite{LR} employed a similar recursion for the moments of the first passage time from level $c$ to $0$, which led them to another derivation of  (\ref{D2}),
the third moment formula
$$
{\mathbb E}[D_0^3]= \frac{6e^{\rho}}{\theta\mu^2}\left[   e^{2\rho} \left(\sum_{j=1}^\infty \frac{\pi_j}{j}   \right)^2 +e^{\rho}\sum_{j=1}^\infty \frac{\pi_j}{j^2}    \right],
$$
and a similar more complex formula for the fourth moment.

Knessl and Young \cite{KY} (p. 217) give a representation of the density of $D_c$ as a series $\sum_{i=1}^\infty c_i \exp(-z_i t/\mu)$
where $z_i$'s are the (positive) roots of $M(-z, c+1-z,\rho)=0$. For example, for $c=0, \theta=\mu=1$ this gives the leading exponential term of the order $\exp(-z_1 t)$ with $z_1=0.450\ldots$,
as compared with $\exp(-t)$ tail of the service time.

For $X(0)=c+1$ the first passage time to $0$ can be represented as  $\sum_{j=0}^c D_j$
with independent $D_j$.
From (\ref{RTA}) and tail asymptotics of the Poisson distribution (cf \cite{Glynn}, Corollary 1 (ii))  for large $c$ we have 
$${\mathbb E}[D_c]\sim \frac{\pi_{c+1}}{\theta\pi_c}=\frac{1}{\mu(c+1)}$$
which gives
$${\mathbb E}\left[\sum_{j=0}^c D_j\right]\sim \frac{\log c} {\mu}, ~~~c\to\infty.$$
Robert \cite{R} (Proposition 6.8) employs the Laplace transform to show that this asymptotics also holds in probability.

See \cite{KY, R} and references therein for asymptotic results in the heavy traffic limit $\rho\to\infty$.

%Knessl and Young \cite{KY} give approximations to the density and the mean of  $D_c$ under various  asymptotic regimes for $\rho$ and $c$.
%Large $\rho$ results are referred to as heavy traffic limits.
%See Section 6.5 of \cite{R} for a  functional limit for $X(\cdot)$, properly scaled and centered, in the asymptotic regime $X(0)\sim\rho x, ~\rho\to\infty$.

\section{Multivariate excursions from the zero state}

A path segment of  $\bc_k(\cdot)$ that  starts  with $(1,0,\ldots,0)$ and terminates upon reaching the zero state $(0,\ldots,0)$ is analogous
to a  busy period of a tandem of  M/M/$\infty$ queues with $k$ phases $S_1,\ldots,S_k$. 
To study the basic characteristics of such multivariate excursion it is enough to follow
 the total 
$$Y(t):=C_1(t)+\cdots+C_k(t),$$
which itself is the occupancy process 
of a single-phase M/G/$\infty$ queue  with Poisson arrival  rate $\theta$, 
and the generic service time $\sigma$ having distribution function
\begin{equation}\label{GG}
{\mathbb P}[\sigma \leq t]=(1-e^{-t})^k.
\end{equation}
Indeed, $\sigma$ is distributed like a sum of exponential variables $\xi_1/1+\cdots+\xi_k/k$, as in (\ref{OrderSt}),
which is  the time that a cycle needs to pass  through $S_1,\dots,S_k$.

The definition of  the busy period (excursion above $0$) for M/G/$\infty$  requires some care, because the process is not Markovian and
the periods spent by present tasks in service cannot be ignored hence must be included in description of  the state \cite{Rolski}. With this in mind, the excursion starting at time $t_0$ is defined under the assumption that
$Y(t_0-)=0$ and $Y(t_0)=1$. In this section we denote by $D_0, H_0, A_0, \Delta_0$ the duration, height, overflow and the number of new arrivals during the busy period of $Y(\cdot)$.

Let $\rho:=\theta\, {\mathbb E}[\sigma]=\theta\, h_k$. The steady state distribution is Poisson$(\rho)$, and formulas (\ref{Wald}), (\ref{RTA}) and (\ref{Edelta}) with $c=0$ apply without change. In particular, 
$${\mathbb E}[D_0]=\frac{e^\rho -1}{\theta}.$$

Recall a concept from the renewal theory.
For a nonnegative integrable random variable $\eta$, representing the generic inter-arrival time,
 the   variable  $\eta^*$ with the integrated tail distribution
$${\mathbb P}[\eta^*\leq t]= \frac{1}{{\mathbb E}[\eta]}\int_0^t {\mathbb P}[\eta>x]{\rm d}x.$$
 appears as the stationary residual lifetime. 
In terms of their Laplace transforms, the relationship between $\eta$ and $\eta^*$ is
\begin{equation}\label{LT*}
{\mathbb E} [\exp(-z \eta^*)]= \frac{1-{\mathbb E}[\exp(-z\eta)]}{z\,{\mathbb E}[\eta]}.
\end{equation}
We shall use this connection of $\sigma$ and $D_0$ to their associated variables $\sigma^*$ and $D_0^*$.

The transient state distribution
$
\pi_c(t):={\mathbb P}[Y(t)=c\,|\,Y(0)=0]
$
is Poisson with mean $\rho\,{\mathbb P}[\sigma^*\leq t]$, hence in particular
$$\pi_0(t)=\exp\{-\rho\, {\mathbb P}[\sigma^*\leq t]\}.$$
From  (\ref{GG}) one finds readily that  $\sigma$ and $\sigma^*$ both  have exponential tails: as  $t\to\infty$ 
\begin{equation}\label{sig-sig}
{\mathbb P}[\sigma>t]\sim k e^{-t}, ~~~{\mathbb P}[\sigma^*>t]\sim \frac{k}{h_k}\,e^{-t},
\end{equation}
which implies that
$$\pi_0(t)-\pi_0=\exp\{-\rho {\mathbb P}[\sigma^*\leq t]\} -e^{-\rho}= e^{-\rho} \big(\exp\{\rho \,{\mathbb P}[\sigma^*> t]\}-1\big)\sim
\theta k e^{-\rho}\, e^{-t},
 $$
where $\pi_0=\lim_{t\to\infty}\pi_0(t)=e^{-\rho}$.

The Laplace transform of the duration is given by
\begin{equation}\label{Ta}
{\mathbb E}[\exp(-z D_0)]= 1+\frac{z}{\theta}- \frac{z}{\theta\,L(z)},
\end{equation}
where
\begin{eqnarray}\label{LI}
L(z)&= &1+\int_0^\infty e^{-zt}   \pi_0'(t) {\rm d}t.
\end{eqnarray}
Equation (\ref{Ta}) is a version of  the
Tak\`acs  formula (\cite{T}, Equation (2) on p. 210)  for the Laplace transform of the time between beginnings of two successive busy periods. 
In \cite{T} and subsequent work
(e.g. Equation (5) in \cite{Stadje}, Equation (4.6) in \cite{LS}) the authors use $L(z)/z$, which is the Laplace transform of $\pi_0(t)$. 
The form (\ref{LI}) is better suitable for our purpose since  $L(z)$ is holomorphic in a larger halfplane
$\Re z>-1$, as dictated by the asymptotics $|\pi_0'(t)|=O(e^{-t})$ for $t\to\infty$.

%where
%\begin{eqnarray}\label{Ldef}
%L(z)&:=&z\int_0^\infty e^{-zt}   \pi_0(t) {\rm d}t.
%\end{eqnarray}
The second moment of the duration was derived from  (\ref{LI}) in Liu and Shi \cite{LS} (Equation (4.13)) as
$${\mathbb E}[D_0^2]= \frac{2}{\theta \pi_0^2}\int_0^\infty (\pi_0(t)-\pi_0){\rm d}t\,.$$
For $k=1$ this has a series representation (\ref{D2})
% (see \cite{Mandjes}, Section 5.2.3),
but for $k>1$ there does not seem to exist a simple analogue.
To compare the numerics, for $\theta=1$
we get  ${\rm Var}[D_0]$ about $12.7921$ for $k=2$ and  about    $4.2123$ for $k=1$.
The joint Laplace transform of $D_0$ and $\Delta_0$ is found in Shanbhag \cite{Shanbhag} (Theorem 2).

We turn next to the counterpart of (\ref{sig-sig}) for the duration of excursion above zero. To that end, designate 
$$F(t):={\mathbb P}[D_0\leq t], ~~~f^*(t):={\mathbb P}[D_0^*\in {\rm d}t]/{\rm d}t,$$
which are the distribution function of $D_0$ and the density function of $D_0^*$, respectively. These are related via 
\begin{equation}\label{fF}
f^*(t)=\frac{1-F(t)}{{\mathbb E}[D_0]}.
\end{equation}

The function $L(z)$ increases from $-\infty$ to $e^{-\rho}$ as $z$ runs from $-1$ to $0$, therefore there exists a unique $\beta\in(0,1)$ satisfying $L(-\beta)=0$
For instance,   $\beta =0.2734\ldots$ if $\theta=1, k=2$.

\begin{theorem}  As $t\to\infty$, it holds that
\begin{equation}\label{exp-as}
1-F(t)\sim \alpha e^{-\beta},
\end{equation}
where
$$\alpha:=  -\left(\theta \int_0^\infty e^{\beta t} \,t \,\pi_0' (t)\,{\rm d}t \right)^{-1}.$$

\end{theorem}

\proof We shall apply a result from the renewal theory. Following
Makowski \cite{Makowski},   the Tak\`acs formula (\ref{Ta}) amounts  to the representation of $D_0^*$ as a geometric sum
$$D_0^*\stackrel{d}{=}  \sum_{j=1}^Q U_j\,,$$
where all variables  involved are independent, $Q$ has the geometric distribution
$${\mathbb P}[Q=j]=\pi_0 (1-\pi_0)^{j-1}, ~~~j=1,2,\ldots$$
%is a geometric distribution with $p=e^{-\rho}$, 
and the $U_j$'s are i.i.d. with  density
$$u(t):= \frac{-\pi^\prime_0(t)}{1-\pi_0}.$$

%$${\mathbb P}[U_j\leq t]=\frac{ 1- \pi_0(t)   }{1-\pi_0}.$$
Conditioning on $U_1$ we arrive at  the improper renewal equation
$$f^*(t) = \pi_0 u(t) +(1-\pi_0)\int_0^t f^*(t-s) u(t) {\rm d}t,$$
with substochastic density $(1-\pi_0)u(t)$.
Adopting a formula from Resnick \cite{Resnick} (page 258, bottom equation where $z(\infty)=Z(\infty)=0$ should be set due to $\lim_{t\to\infty} \pi_0'(t)=0$)
we have
$$f^*(t)\sim \alpha^* e^{-\beta t},$$
with  $\beta\in (0,1)$ as above solving $L(-\beta)=0$ (cf \cite{Resnick}, Proposition 3.11.1) and

\begin{eqnarray*}
\alpha^* &=&\frac{\pi_0\int_0^\infty e^{\beta x} u(x){\rm d} x}  {(1-\pi_0)\int_0^\infty x e^{\beta x}u(t) {\rm d}t}\\
&=&\frac{\alpha\theta}{e^{\rho}-1}.
\end{eqnarray*}
The assertion now follows by the virtue of (\ref{fF}).
\endpf

An alternative approach 
%to the exponential tail asymptotics
 is the following.
Using (\ref{LT*}) we have 
\begin{equation}\label{LTf}
{\mathbb E}[\exp({-z D_0^*})]=\int_0^\infty e^{-zt} f^*(t){\rm d}t=(e^\rho -1)^{-1} \left( \frac{1}{L(z)}-1\right).
\end{equation}
From this the exponential tail asymptotics can be concluded by  singularity analysis of the Laplace transform.
Indeed, with $\Re z$ fixed, $|L(z)-1|$ is maximised for $\Im z=0$, hence and by monotonicity $L(z)\neq 0$ if $\Re z>-\beta$. 
On the other hand, 
by a property of the Laplace transform
$|L(z)-1|\to 0$ as $|z|\to \infty$ uniformly in $\Re  z>-\beta-\varepsilon$
(cf \cite{Doetsch}, Theorem 23.6). Thus for $\varepsilon>0$ sufficiently small,  $L(z)$ has no zeros in this halfplane other than $-\beta$, hence 
the only singularity of (\ref{LTf}) in the halfplane is   a simple pole  at $-\beta$, with residue readily identified with $\alpha^*$.
From this  (\ref{exp-as}) follows by writing $f^*(t)$ in the form of the inverse Laplace transform of (\ref{LTf}), then moving the contour of integration to $\Re z=-\beta-\varepsilon$,
see \cite{Doetsch} (Section  35) for this classic technique.

Note that $\alpha$ is the residue of (\ref{Ta}) at pole  $-\beta$, but using (\ref{Ta}) directly to justify the tail asymptotics (of the density of $D_0)$  looks more difficult due to the factor $z$.

\section{The embedded   tagged cycle process}

Suppose element $n$ starts a new cycle of the CRP permutation, with some number $L_1^{(n)}$  of singletons already present in $\Pi^{(n-1)}$,
that is $L_1^{(n)}=C_1^{(n-1)}$.
Let $L_2^{(n)}$ be the number of doubletons present 
immediately before this cycle moves to $S_2$, etc.
Intuitively, $L_k^{(n)}$  is what an observer moving with the tagged cycle spots in $S_k$ when entering the phase. 
As $n\to\infty$ the distribution of  $L_1^{(n)}, L_2^{(n)},\ldots$ converges to a limit which 
has Poisson marginals as 
in (\ref{CZn}) but they are not independent.
It seems hard to capture features of the limit multivariate distribution without turning to the  embedding of CRP in continuous time.
Fortunately,   a major work has been done by the queueing theorists.

To set a general scene,
consider a  tandem  of M/M/$\infty$ queues with arrival        rate $\theta$   and sojourn parameter $\mu_k$ for phase $S_k$, and let $\rho_k:=\theta/\mu_k$.
 Assuming the system in steady state and that there is a tagged arrival at time $0$,           
let $ L_k$ be
the occupancy of          $S_k$  immediately before the   tagged item enters $S_k$.
The following result was obtained by Vainstein and Kreinin \cite{VK1} and extended by Boxma \cite{Boxma} to  tandems of M/G/$\infty$ queues with arbitrary sojourn times.
As above,  $\xi_1,\xi_2,\ldots$ denote independent unit exponential random variables.

\begin{theorem} The joint p.g.f.    of $L_j$ and $L_k$ for $1\leq j<k$ is
\begin{eqnarray}\label{LL}
{\mathbb E}[x^{L_j} y^{L_k}]=\exp\{\rho_j(x-1)+\rho_k  (y-1)\}\int_0^\infty \exp\{\rho_j(x-1)(y-1)\varphi_{jk}(t)\}{\rm d}\psi_{jk}(t),
\end{eqnarray}
where
$$ \varphi_{jk}(t) := {\mathbb P}\left[ \frac{\xi_j}{\mu_j}+\cdots+     \frac{\xi_{k-1}}{\mu_{k-1}}<t<\frac{\xi_j}{\mu_j}+\cdots+     \frac{\xi_{k}}{\mu_{k}} \right],~~~
 \psi_{jk}(t) := {\mathbb P}\left[ \frac{\xi_j}{\mu_j}+\cdots+     \frac{\xi_{k-1}}{\mu_{k-1}}<t\right].$$

\end{theorem}
\proof 
By the steady-state assumption, the flow from $S_{j-1}$ to $S_j$ is Poisson, hence
we will not lose generality by considering the case $j=1$ only.
We adapt the more general argument from \cite{Boxma} (Theorems 2.2 and 3.1) to the  M/M/$\infty$ tandem.
We have $L_1\stackrel{d}{=}{\rm Poiss}(\rho_1)$, and the time, say  $T$, that the tagged arrival to $S_1$ needs to reach $S_k$
has distribution function ${\mathbb P}[T\leq t]=\psi_{1k}(t)$. For $1\leq i\leq k,$ let $N_i$ be the number of items  in $S_k$ at time $T$
that were in $S_i$ at time $0$,
and let $N_0$ be the number of items in $S_k$
at time $T$ that were not yet present in the system at time $0$. Clearly, $L_k=N_0+N_1+\cdots+N_k$.
Given the tagged item finds $L_1=\ell_1$ and needs time $T=t$,  the variables $N_0,\ldots,N_k$ are conditionally independent,
$$N_0\stackrel{d}{=}{\rm Poiss}(\theta p_0),~~ N_1\stackrel{d}{=}{\rm Bin}(\ell_1, p_1){\rm ~~and~~} N_i\stackrel{d}{=}{\rm Poiss}(\rho_i p_i),~2\leq i\leq k.$$
Here,
\begin{eqnarray*}
p_0&=&\int_0^t \varphi_{1k}(x){\rm d}x,
\end{eqnarray*}
and $p_i$ for $1\leq i\leq k$ is the   probability that the generic item from $S_i$ is located in $S_k$ over time $t$, i.e.
\begin{eqnarray*}
p_i& =&{\mathbb P}\left[ \sum_{m=i}^{k-1}     \frac{\xi_m}{\mu_m}<t<\sum_{m=i}^{k}     \frac{\xi_m}{\mu_m} \right].
\end{eqnarray*}
Note that $p_1=\varphi_{1k}(t)$.
The steady-state balance equation for the mean content of $S_k$ is
$$\theta p_0+\rho_1 p_1+\cdots+\rho_k p_k=\rho_k,$$
which allows us to write 
$\theta p_0+\rho_2 p_2+\cdots+\rho_k p_k=\rho_k-\rho_1 p_1,$
and together with the above conclude that the  conditional distribution of $L_k$ is the convolution
$${\rm Bin}(\ell_1,p_1)*{\rm Poiss}(\rho_k-\rho_1 p_1),$$
whence
$${\mathbb E}[x^{L_1}y^{L_k}\,|\,L_1=\ell_1, T=t]= x^{\ell_1}\{1-\varphi(t)+\varphi(t) y\}^{\ell_1} \exp\{( \rho_k-\rho_1\varphi(t) ) (y-1)\}. $$
The result now  follows by integrating out $\ell_1$ and $t$.                                                                                                                                                                                                                                       
\endpf

From (\ref{LL})   follows that (or see \cite{Boxma}, Equation (3.4))
\begin{equation}\label{cov}
{\rm cov} (L_j,L_k)=\rho_j \int_0^\infty \varphi_{jk}(t){\rm d} \psi_{jk}(t), ~~~{\rm corr}(L_j,L_k)=  \sqrt{\frac{\mu_k}{\mu_j}} \int_0^\infty \varphi_{jk}(t){\rm d} \psi_{jk}(t).
\end{equation}
Vainshtein and Kreinin \cite{VK2} observed that 
\begin{equation}\label{Lag}
{\rm corr}(L_j, L_k)= \frac{1}{2\sqrt{\mu_j\mu_k} }\,{\cal  L}(0),
\end{equation}
where ${\cal L}(\cdot)$ is the Lagrange polynomial interpolating the square root function from  the data set $(\mu_j^2,\mu_j),\ldots, (\mu_k^2,\mu_k)$.
Remarkably, the correlation coefficient does not depend on $\theta$.

 The case relevant to permutations
 \begin{equation}\label{muk}
\mu_k=k
\end{equation}
will be worked out in the rest of this section. 
Using (\ref{Lag}) Vainstein and Kreinin 
(\cite{VK2}, Equation (15)) evaluated (\ref{cov}) 
as 
$${\rm corr}(L_1,L_k)=\frac{\sqrt{k}}{2(2k-1)}.$$
 We take a different approach which works smoothly for all $j$ but is limited to  (\ref{muk})  (or constant multiples of (\ref{muk})).
Let $\xi_{j:k},~1\leq j\leq k$, denote the $j$th {\it maximal} order statistic among  the first $k$ exponential variables.
Using Renyi's representation we have the identities
\begin{eqnarray*}
\psi_{jk}(t)&=&{\mathbb P}[\xi_{j:k-1}<t],\\
\varphi_{jk}(t)&=&{\mathbb P}[\xi_{j:k-1}<t]-{\mathbb P}[\xi_{j:k}<t]=
{\mathbb P}[\xi_{j:k-1}<t<\xi_{j-1:k-1}]\,{\mathbb P}[\xi_k>t],
\end{eqnarray*}
where the last equality follows from the events coincidence 
$$\{  \xi_{j:k-1}<t, \, \xi_{j:k}\geq t \} =    \{\xi_{j:k-1}<t\leq  \xi_{j-1:k-1},\, \xi_k\geq t \}.$$
To express (\ref{cov}) via a beta integral we    pass to the uniform  order statistics, thus obtaining
\begin{eqnarray*}{\varphi_{jk}}(-\log(1-x))&=& {k-1\choose j-1} x^{k-j} (1-x)^{j},\\
{\rm d}{\psi_{jk}}(-\log(1-x))&=& (k-1){k-2\choose j-1}x^{k-j-1}  (1-x)^{j-1}   \,      {\rm d}x\,,
\end{eqnarray*}
whence     (\ref{cov}) for  rates (\ref{muk}) becomes

\begin{eqnarray*}
{\rm cov}(L_j,L_k)&=&\,~~~~~\theta{k-1\choose j}{k-1\choose j-1}\frac{(2k-2j-1)!(2j-1)!}{(2k-1)!}\,,\\
{\rm corr}(L_j,L_k)&=&\sqrt{jk}{k-1\choose j}{k-1\choose j-1}\frac{(2k-2j-1)!(2j-1)!}{(2k-1)!}\,.
\end{eqnarray*}
Interestingly, the covariance has some symmetry, ${\rm cov}(L_j,L_k)={\rm cov}(L_{k-j},L_k)$.

Since $\varphi_{jk}(t)=\psi_{j,k}(t)-\psi_{j,k+1}(t)$, (\ref{cov}) implies an estimate
$${\rm corr}(L_j,L_k)<\frac{1}{2}\sqrt{\frac{j}{k}}\,,$$
which gives the correct decay order $k^{-1/2}$ of the correlation as $k\to\infty$ and $j$ is fixed.

\section{A functional limit for the small cycle counts}

Finally, we argue that $(\bc_k(t),~t\geq 0)$ appears as a weak limit of  $(\bc_k^{(n)},~n\geq 0)$ by the virtue of a nonlinear time-scale change.
To that end, we interpolate the discrete time Markov chain to a piecewise constant jump process with real time  parameter.

\begin{theorem} Let $\bc_k(\cdot)$ start at time $0$ in some random state $\bc_k(0)$,
and let $\boldsymbol{C}_k^{(\cdot)}$ start at time $\nu$ in some   random state $\bc_k^{(\nu)}$.
If $\bc_k^{(\nu)}$ converges in distribution to $\bc_k(0)$ as $\nu\to\infty$ then also
$$(\boldsymbol{C}_k^{( \nu e^t)}, ~t\geq 0)\Rightarrow (\bc_k(t),  ~t\geq 0)\,,$$
where $\Rightarrow$ means weak convergence  in the Skorohod space $D[0,\infty)$.
\end{theorem}

To ease notation let $X_\nu(t)=\boldsymbol{C}_k^{( \nu e^t)}$.
Since the   state space is discrete,   the assertion can be reduced    to the  case when the initial states   are fixed and identical, that is $\boldsymbol{C}_k(0)=X_\nu(0)={\boldsymbol c}(0).$
The embedded    jump    chains  have the same   transition  probabilities, hence    it is possible to couple the processes   in such a way that they pass the same random sequence of states.  
Appealing to \cite{Xia} (Lemma 2.12) shows
that is suffices to verify that the sequence of consecutive sojourn times of $X_\nu(t)$, seen as a random element of ${\mathbb R}_+^\infty$,  converges  in  distribution
to the sequence of sojourn times of $\boldsymbol{C}_k(t)$.
Given a path ${\boldsymbol c}(0), {\boldsymbol c}(1),\ldots$ of the jump chain, the sojourn times of $\boldsymbol{C}_k{(\cdot)}$ are independent  exponential variables, with rates         
\begin{equation}\label{rates}
r=\theta+\sum_{i=1}^k i c_i
\end{equation}
   depending on  ${\boldsymbol c}\in {\mathbb Z}_+^k$.
The next lemma finds the limiting distribution of the sojourn time of $X_\nu(t)$ at 
an arbitrary state ${\boldsymbol c} \in {\mathbb Z}_+^k$.
\begin{lemma} Given $X_\nu(t)=\boldsymbol{c}$ the residual sojourn time in this state converges in distribution to ${\rm Exp}(r)$, as $\nu\to\infty$.
\end{lemma}

\proof Using (\ref{calc}), as $n\to\infty$ we obtain
$$
{\mathbb P}[\bc_k^{(i)}=\boldsymbol{c}\,,n\leq i\leq   m\,|\,\bc_k^{(n)}=\boldsymbol{c}]=
\prod_{j=n}^{m-1} \left(1- \frac{r}{j+\theta}\right)=\left(\frac{n}{m}\right)^r\left(1+O\left(\frac{1}{n}\right)\right)
$$
uniformly in $m\geq n$.
Setting $n=\nu e^t, ~m=\nu e^{t+\delta}$, $\delta>0$,  we conclude that $X_\nu(\cdot)$ spends in $\boldsymbol{c}$ some  time 
exceding 
 $\delta$  with probability $e^{-r\delta}+O(\nu^{-1})$,  hence the limit distribution is exponential as stated.
\endpf

Let $r(0), r(1),\ldots$ be the rates for $\boldsymbol{c}(0), \boldsymbol{c}(1), \ldots$ defined by (\ref{rates}),
and let $V_\nu(0), V_\nu(1),\ldots$ be the sojourn times that $X_\nu(\cdot)$ spends  in these states.
By the lemma, $V_\nu(0)$ converges in distribution to ${\rm Exp}(r(0))$. By the strong Markov property and because the estimate $O(\nu^{-1})$ in the proof of the lemma is uniform in $t$,
the conditional distribution of  $V_\nu(1)$ given $V_\nu(0)$ converges to ${\rm Exp}(r(1))$. But then we also have the joint convergence of $(V_\nu(0), V_\nu(1))$,
as follows from \cite{Sethuraman}  (Theorem 2). Continuing by induction, we obtain the joint convergence of the sojourn times $V_\nu(0), V_\nu(1),\ldots$ to the counterpart 
sequence of  sojourn times of $\bc_k(\cdot)$ and the weak convergence of $X_\nu(\cdot)$ follows.
\endpf

\end{document}